\documentclass{ws-rv9x6}  

\def\picill#1by#2(#3)
{\vbox to #2
{\hrule width #1 height 0pt depth 0pt
\vfill\epsffile{#3}}}

\begin{document}

\setcounter{chapter}{0}

\chapter{Minimal Flat Knotted Ribbons}

\markboth{L. Kauffman}{Minimal Flat Knotted Ribbons}

\author{Louis H. Kauffman}

\address{Department of Mathematics, Statistics and  Computer Science
(m/c 249)\\851 South Morgan Street\\University of Illinois at Chicago\\
Chicago, Illinois 60607-7045\\
E-mail: kauffman@uic.edu}

\begin{abstract}
This paper gives mathematical models for flat knotted ribbons, and makes specific conjectures for the least length of 
ribbon (for a given width) needed to tie the trefoil knot and the figure eight knot. 
 
\end{abstract}

\section{Introduction}     
We study framed knots, configured as knotted ribbons, where a ribbon is a space homeomorphic to a rectangle.
Think of that rectangle as a long strip of paper, and imagine tying a knot in the strip and then pulling it gently tight
so that it becomes a flat folded knotted strip. See Figure 3 for an illustration of this process for the trefoil knot.
The folds in the strip should all be of the type shown in Figure 2. The question about such knots is this: What is the least 
length of strip (for a given width) needed to make a flattened version of the knot? 
\bigbreak 

This question, about length to width ratio for flat knots, is analogous to a question about length to radial width (so-called {\em
thickness}) for knotted tubes embedded in three dimensional space. The question for knotted tubes has a good-sized recent
literature \cite{BS,CKS,Rawdon,Ideal}, but to this date there is no specific conjecture for the exact thickness of even the trefoil
knot. In the present  paper, we make specific conjectures for minimal length to width ratio for the trefoil and figure eight as knotted
flat ribbons. It is clear from this initial exploration of knotted ribbons that there is much to think about and very interesting
relationships between geometry and topology in this domain.
\bigbreak

The paper is organized as follows. Section 2 gives our mathematical model for flat knotted ribbons, in terms of immersed piecewise linear
ribbons in the plane. We emphasize two approaches to the length to width ratio, one for closed knots, and one for truncated knots.
In Section 3 we analyze the geometry of the folded trefoil and conjecture its minimal ratio for the truncated trefoil knot. In Section 4,
we analyze the geometry of the folded figure eight knot and conjecture its minimal ratio. Section 5 is a discusssion of possible
relationships with fully three dimensional issues.
\bigbreak

\section{A Mathematical Model for Flat Knotted Ribbons}
In Figure 1 we depict two knotted ribbons ($A$ and $B$) embedded in three dimensional space. The projection of ribbon $A$ to the
plane is {\em not} the sort of flattened knot that we are considering here. The twist in the embedding flattens to a singular
projection that we  wish to avoid. The projection of ribbon $B$ to the plane is suitable for our purposes.
\bigbreak

$$ \picill3inby1.5in(knottedband) $$

\begin{center}
{\bf Figure 1.  Knotted Ribbons}
\end{center}

$$ \picill3inby3.5in(MirrorsFolds) $$

\begin{center}
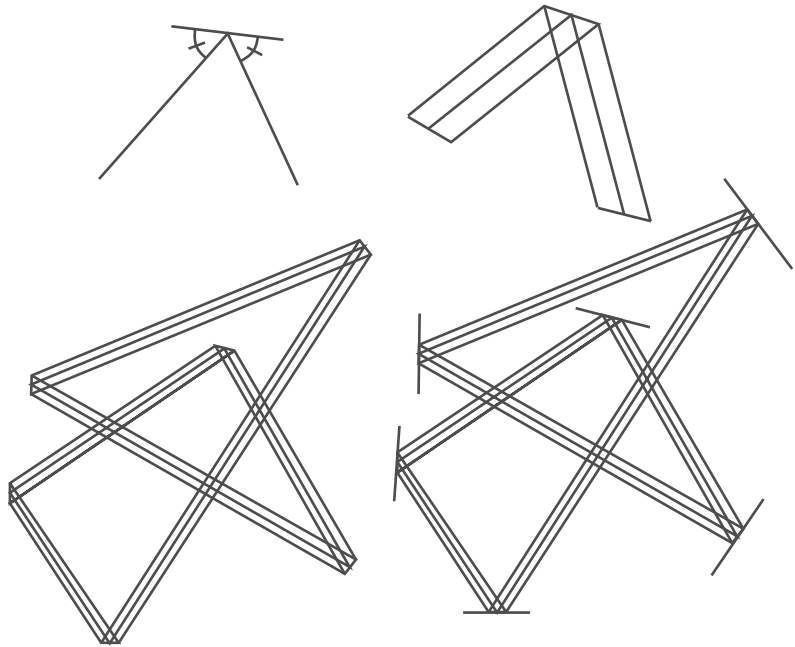
{\bf Figure 2.  Mirrors and Folds}
\end{center}

In Figure 2 we provide a mathematical model for the 
sort of folded strip that we have in mind.  In the mathematical model, there is a piecewise linear curve in the plane with 
transverse self-intersections. This curve is the projection of the center-line of the strip to the plane after all folds have 
been made. Each fold is modeled as follows: Consider an angle in the piecewise linear embedding. Place a {\em mirror segment} at each
angle so that it makes equal angles to each of the parts of the embedding that form the angle. The two parts of the embedding at the 
angle now have the appearance of rays being reflected by the mirror segment. Then the image of the folded ribbon in the plane is {\em by
definition} a bundle of parallel rays reflected by the mirror segments. See Figure 2 for an illustration of this geometry.
It is assumed that separate
straight line bundles intersect either along mirrors, or transversely in the pattern of the self-crossings of the core diagram. A {\em
flat ribbon immersion } is such an immersion of a ribbon into the plane. A {\em flat knotted ribbon} is a choice of weaving that
overlies a flat ribbon immersion. In such a weaving, the width of any straight interval of parallel rays must be chosen so that {\em no
bundle of parallel lines can overlap a fold line when it is intended to weave that interval with the intervals meeting at the fold}.  Such
an overlap would be the same as a situation where the bundle of lines overlapped one of the mirrors, and would have to be either entirely
over, or entirely under the mirror fold. Such a choice is obtained by  first selecting an over or under crossing at each
self-intersection of the center-line, and then extending these choices compatibly across the families of reflected rays. In this way, any
piecewise linear knot diagram gives rise to a collection of flat knotted ribbons of varying width. For a given diagram, there will be a
unique largest width knotted ribbon. We would like, within the topological type of a given knot, to determine the least length of ribbon
needed  for a given choice of width.
\bigbreak

In order to discuss the issue of the length of a given knotted ribbon (for a given width), there are two choices. We can assume that the 
knotted ribbon is in the traditional form of a closed piecewise linear loop, and take the length of the core curve. This is the {\em
closed length} of the knotted ribbon. Alternatively, we can choose to cut and truncate the knotted ribbon and measure the length of the
unfolded ribbon that results from the cut. In some cases, there is a clear way to make the truncation. For example, in Figure 3 we
illustrate how a trefoil knot occurs in the form of a woven pentagon. In this form it is natural to truncate the ribbons flush with the
sides of the  pentagon. When there is such a choice, we shall call the resulting length the {\em truncation length} of the knotted
ribbon.

In Figure 3, we illustrate the process of forming a trefoil knot and pulling it to make a pentagonal flat knotted ribbon.
We also show the internal mirror structure of this knot. The edges of a regular pentagon are the mirror segments. A ray of light
enters the pentagon and is reflected three times before exiting. In the diagram, the over and undercrossings of the trefoil knot
have been indicated. Note that the pentagonally folded trefoil knot is seen here as a limiting case of flat banded knots whose mirrors
are segments on the sides of a regular pentagon. There is no way, using this particular geometry, to extend the lengths of the mirrors
beyond the lengths of the pentagonal sides and still maintain the possiblity of an embedding.
\bigbreak 

We also show in Figure 3 how a totally internally self-reflecting pentagon would behave. In this case the alternating over
and undercrossing pattern is actually not realizable by straight lines in three space, but could be effected by a weaving process,
as we have described.
\bigbreak

$$ \picill3inby5in(Mirrors) $$

\begin{center}
{\bf Figure 3.  A Pentagonal Mirror}
\end{center}

\section{The Trefoil Knot, the Pentagon and the Golden Ratio}

We begin with an experiment that one can perform with the pattern of the
trefoil knot and a strip of paper. Tie the strip into a
trefoil and pull it gently tight and fold it so that you
obtain a flat knot. As you pull it tight, a regular pentagon will appear.
\bigbreak

It is intuitively clear that this pentagonal form of the
trefoil knot uses the least length of paper for a given
width of paper (to make a flattened trefoil). At this
writing, we do not have a proof of this statement. It is the purpose of this section to analyze the pentagonal geometry 
and to make a quantitative version of this conjecture.
\bigbreak

$$ \picill4inby5.6in(PullUnfoldTrefoil) $$

\begin{center}
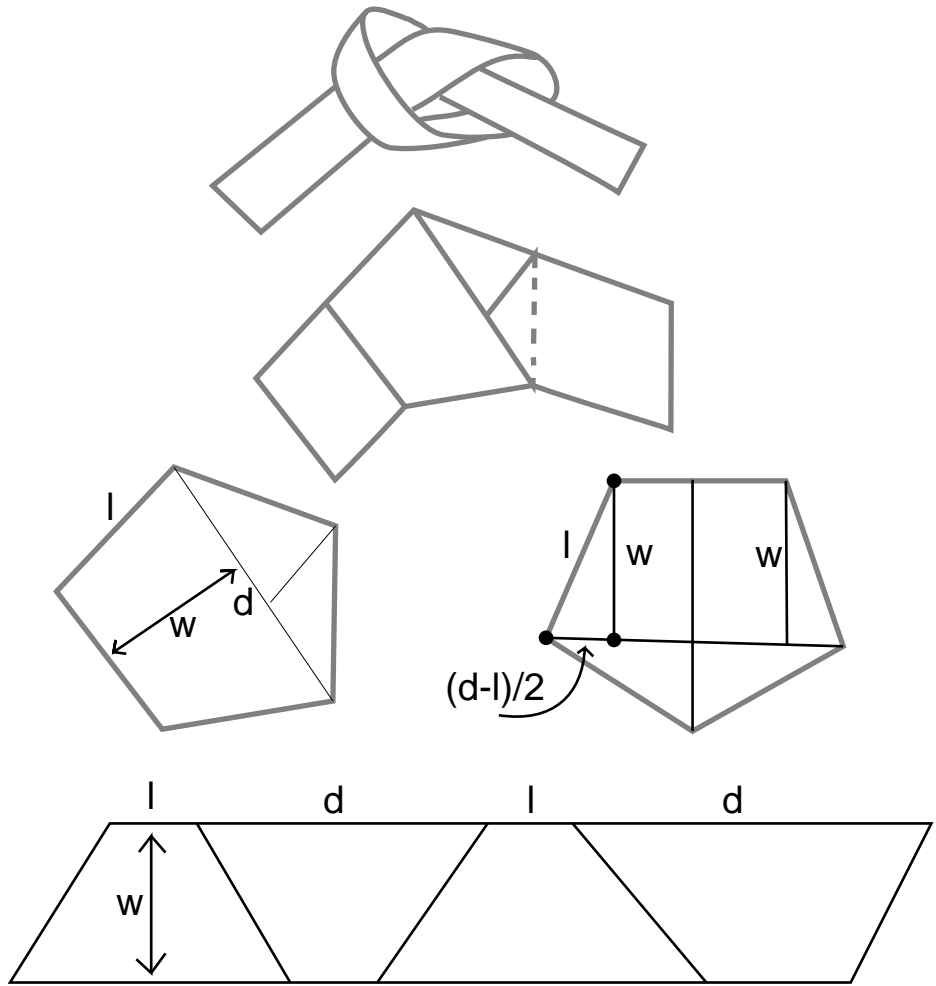
{\bf Figure 4.  Pull and Unfold Trefoil}
\end{center}

We can determine the length to width ratio of
the strip of paper obtained from this pentagonal construction, by first
cutting the knot exactly along the edges of the pentagon, and then unfolding this strip.
See Figure 4.
\bigbreak

The lengths  $l$  and   $d$  that label the  strip in Figure 4 are the
edge length  $l$ and the chord length $d$ of the pentagon.  To
see this, contemplate Figure 4, or better yet, perform the experiment of making the knot, cutting it flush to the 
pentagon and unfolding the strip. The length of the whole strip is  $L = 2(l+d).$
\bigbreak

In Figure 4 there is a blown-up picture of the pentagon with the
edge-length  $l$, the chord $d$ and the width of the strip $W$
labeled. You can see that the
width is indeed the length of a particular perpendicular dropped to the
chord. It is also apparent from this diagram that  
$$W^{2}  + (d-l)^{2}/4 = l^{2}.$$ Thus we
can solve for the width $W$  of the strip in terms of $l$ and $d.$
The result is

$$W = \sqrt{3l^{2} - d^{2} + 2ld}/2.$$
Since  the length of the strip is $L = 2(d+l),$ the ratio of length to width is

$$L/W = 4(d+l)/\sqrt{3l^{2} - d^{2} + 2ld}$$

$$= 4((d/l) +1)/\sqrt{3 - (d/l)^{2} + 2(d/l)}.$$
Let $$\phi = d/l.$$ Then we have

$$L/W = 4(\phi+1)/\sqrt{3 - \phi^{2} + 2\phi}.$$
\bigbreak

This is a formula for the length to width ratio of the
flattened (truncated) pentagonal trefoil knot.  Now we use a
classical fact about the pentagon.
\bigbreak

\noindent {\bf Fact.}  The ratio $\phi = d/l$ of the chord of a regular pentagon
to its edge length is equal to the golden ratio.  That is, 
$\phi = (1 + \sqrt{5})/2$ and $\phi^{2} = \phi + 1.$
To see this fact, contemplate Figure 5.
\bigbreak

$$ \picill3inby4in(Pent) $$

\begin{center}
{\bf Figure 5.  The Pentagon}
\end{center}

The small pentagon has its edge  $l$ and chord  $d$ labeled. We
embed the small pentagon in the larger one and observe via a
parallelogram and by similar triangles that $d/l = (d+l)/d$. 
Thus, with $\phi = d/l,$ we have $d/l = 1 + l/d$
whence $\phi = 1 + 1/\phi.$ This is sufficient to show that  $\phi$ is the
golden ratio. Note how the golden ratio appears here through
the way that a pentagon embeds in a pentagon, a pentagonal
self-reference.
\bigbreak

We conclude that

$$L/W = 4(\phi+1)/\sqrt{3 - \phi^{2} + 2\phi} = 4(\phi+1)/\sqrt{2 + \phi}$$
Thus
$$L/W = 4(\phi+1)/\sqrt{2+\phi} = 4/\sqrt{7-4\phi}.$$
(We thank Robert Gray \cite{Gray} for observing this last equality.)

\noindent This shows how the golden ratio appears crucially in this
fundamental parameter associated with the trefoil knot.

Numerically, we find that

$$L/W =5.5055276818846941528288383276435507 \cdots$$

\noindent {\em We conjecture that $$L/W = = 4/\sqrt{7-4\phi} = 4/\sqrt{5-2\sqrt{5}}$$
is indeed the minimum length to width ratio for a flattened truncated
trefoil knot.}
\bigbreak

\noindent {\bf Remark.} It is well-known in recreational mathematics
 that the trefoil ribbon produces a pentagon \cite{Schneider}. To our knowledge, this is the first time this construction has been analyzed for its
length to width ratio, and it is the first time that the construction has been placed in the context of the search for 
minimal ribbon length.
\bigbreak

\section{Figure Eight Knot and a Hexagon}

$$ \picill4inby2.2in(FigureEight) $$

\begin{center}
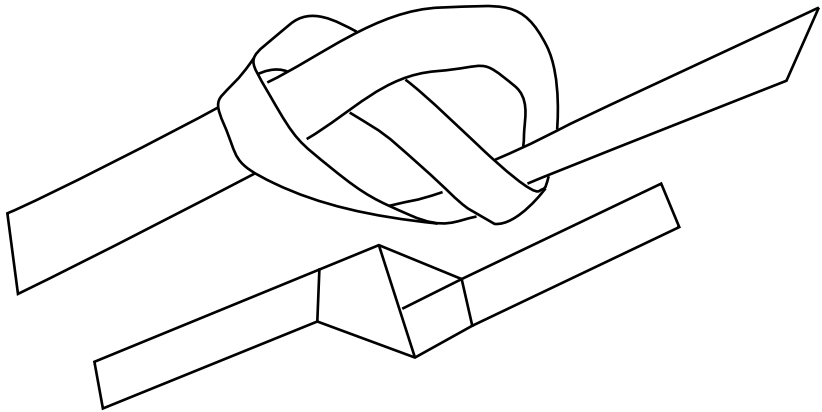
{\bf Figure 6.  Figure Eight Knot and its Hexagon}
\end{center}

In this section we consider the analogous question for the figure eight knot.
Figure 6 illustrates the experiment for the figure eight knot that leads to a flat folded hexagon.
As the figure shows, the hexagon is not regular. It has two opposite parallel sides of the same length, and four sides that are equal in
length to one another. 
\bigbreak

$$ \picill4inby4.5in(EightHex) $$

\begin{center}
{\bf Figure 7.  Geometry of Folded Figure Eight Hexagon}
\end{center}

Consider Figure 7, which shows the hexagon formed by the tightened ribbon and the center line of the
ribbon (appearing as five line segments with four folds), together with a sixth line segment connecting the starting point ($ A $) and ending point ($ H $) of
this folded center line.  
Place the figure on a coordinate plane so that the point $ O = (0,0) $, $ A = (a, b) $, and $ B = (a+c, 0) $.  
Since $ A $ is the midpoint of $ BC $, we must have $ C = (a-c, 2b) $ and $ D = (0, 2b) $.  Finally, since lines $ BC $ and $ DG $ are parallel, 
\[
 G = B + (D - C) = ((a+c) - (a-c), (2b) - (2b)) = (2c, 0) .
\]
On the other hand, $ \angle AGB = \angle ABG $ because $ AG $ is parallel to $ DH $, so $ G = (a-c, 0) $ 
implying that $ a = 3c $.  Thus $ B = (\frac{4}{3} a, 0) $, $ C = (\frac{2}{3} a, 2b) $, and $ G = (\frac{2}{3} a, 0) $.
\bigbreak

\noindent Now the angle of incidence at $ A $ equals the angle of reflection, so $ \angle BAG = \angle CAF $.  
Furthermore,  $ FG $ and $ AB $ are parallel, so $ \angle BAG = \angle AGF $ and $ \angle CAF = \angle AFG.$
Hence $ \angle AFG = \angle AGF.$ 
Thus $\triangle{GAF}$ is isosceles, so that $|AF| = |AG|$.
\bigbreak

\noindent Since $ F $ is the intersection point of the lines $ OA $ and $ DG $, it satisfies the equations 
$$ 
y = \frac{b}{a} x \qquad \text{and} \qquad
y = \frac{-3b}{a} x + 2b .
$$
Hence $ F = (\frac{a}{2}, \frac{b}{2}) $.  
Then $ |AF|^2 = |AG|^2 $ implies 

\[
\frac{1}{4} a^2 + \frac{1}{4} b^2 = \frac{1}{9} a^2 + b^2 , 
\] 
so that 
\[ b = \sqrt{\frac{5}{27}} a .\]

\noindent The center line of the ribbon starts at $ (a, b) $ and proceeds to $ (0, -2b), (-a, b), (a, -b) $, and $ (0, 2b) $, before finishing at $ (-a, -b) $.
 Thus the length of the ribbon is
\[
L = 4 \sqrt{a^2 + 9b^2} + \sqrt{4a^2 + 4b^2}
= \frac{32}{3} \sqrt{\frac{2}{3}} \, a 
\]

\noindent In order to find the width $W$ of the ribbon, view Figure 8. Here we observe that $W/2$ is the distance between the parallel lines through $AB$ and 
$FG.$ It follows from our coordinate calculations that $|AB| = |AG|$ since
$$|AB| = |(a, b) - (4a/3, 0)| = |((a, b) - (2a/3, 0)| = |AG|.$$
From this, and the equality of alternate interior angles between parallel lines, it follows that the three angles indicated by $\alpha$
in Figure 9 are indeed equal to one another. That is, $\angle AGB = \angle ABG = \angle BGK.$ It then follows at once that 
$W = (4/3)a \sin(\alpha)$ and that $\alpha = \arctan(3b/a) = \arctan(\sqrt{5/3}).$ Hence
$$W = (4/3)a \sin(\arctan(\sqrt{5/3})) = (4/3)a \sqrt{5/3}/\sqrt{1 + 5/3} = a \sqrt{10}/3.$$
and so
\[
L/W = \frac{32}{3} \sqrt{\frac{2}{3}} \frac{3}{\sqrt{10}}
= \sqrt{\frac{2}{3}} \frac{32}{\sqrt{10}}
= \frac{32}{\sqrt{15}}
\approx 8.26236447190916 \cdots .
\]
\bigbreak

$$ \picill5inby3.7in(EightHex1) $$

\begin{center}
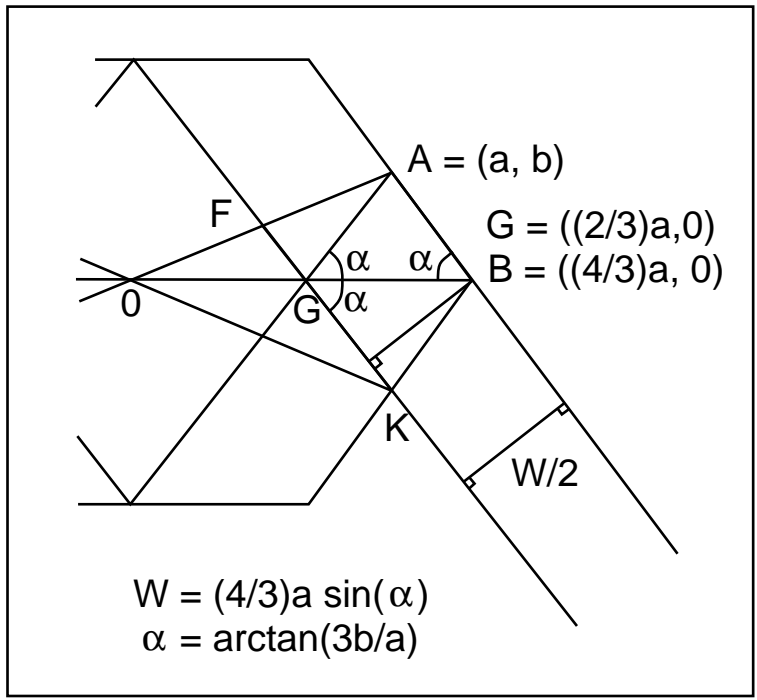
{\bf Figure 8.  Finding the Width W}
\end{center}

\noindent As expected, the figure eight knot is longer than the trefoil knot. We conjecture that $\frac{32}{\sqrt{15}}$ is the minimal (truncation)
length to width ratio for the figure eight knot. 
\bigbreak

Even at the level of the figure eight knot there are other ways to make the folds. But experiments so far reveal that the $L/W$ ratio
is larger for these other constructions. An example of an alternative fold is shown in Figure 9. We leave it to the reader to verify 
that this folding geometry leads to more length than the semi-regular hexagon described above.

$$ \picill3inby1.5in(AltEight) $$

\begin{center}
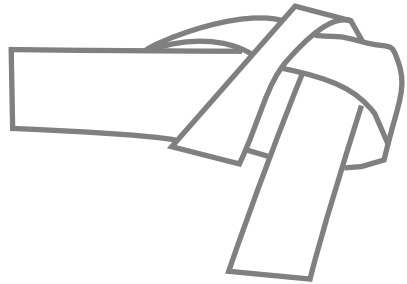
{\bf Figure 9.  An alternate, longer, flat figure eight knot}
\end{center}

\section{Discussion}
Some workers have thought about the relationship between thickness for flat knots and thickness for knotted tubes.
For example Kusner \cite{Kusner} informed the author of this paper that 
one can get an upper bound of the form 

$$ Ropelength(K) \le c{_1}Ribbonlength(K) + c{_2} Crossing(K) $$
He remarks that it seems that the best one can expect for ribbonlength is
a linear relation

$$c_{3} Crossing(K) \le Ribbonlength(K) \le c_{4} Crossing(K)$$
with minimum crossing number. Here $Ropelength$ and $Ribbonlength$ refer respectively to minimal length to width ratios for
tubular and flat knots, while $Crossing$ refers to the minimal crossing number.
\bigbreak

It is clear that much more work remains to be done in this field. The immediate appearance of geomerty, in relation to flattened knotted
ribbons, is encouraging. The problems of minimization discussed in this paper can be explored both theoretically and via computer models.
There should be important and interesting relationships between the structure of ribbons, the structure of ropes, and geometry in 
three dimensional space. Values such as the minimal ropelength and minimal ribbonlength for knots are fundamental topological quanta
associated with the knots. It is to be expected that such numbers will be related to the geometry, topology and physics of these entities.
\bigbreak

\noindent {\bf Acknowledgment.} Most of this effort was sponsored by the Defense
Advanced Research Projects Agency (DARPA) and Air Force Research Laboratory, Air
Force Materiel Command, USAF, under agreement F30602-01-2-05022.  
The U.S. Government is authorized to reproduce and distribute reprints
for Government purposes notwithstanding any copyright annotations thereon. The
views and conclusions contained herein are those of the authors and should not be
interpreted as necessarily representing the official policies or endorsements,
either expressed or implied, of the Defense Advanced Research Projects Agency,
the Air Force Research Laboratory, or the U.S. Government. (Copyright 2004.) 
It gives the first author great pleasure to acknowledge support from NSF Grant DMS-0245588,
to give thanks to the University of Waterloo and the Perimeter Institute in Waterloo, Canada for their hospitality during
the preparation of this research, and to thank Jorge Alberto Calvo for penetrating remarks about the geometry of the figure eight knot.

\end{document}